\documentclass[11pt,a4paper]{article} %ГришаТамасян
\usepackage[cp866]{inputenc}
\usepackage[english]{babel}
\usepackage{amsmath}
\usepackage{amsfonts}
\usepackage{longtable}
\usepackage{amssymb}
\usepackage{graphics}
\newcommand{\bc}{\begin{center}}
\newcommand{\br}{\begin{right}}
\newcommand{\ec}{\end{center}}
\newcommand{\be}{\begin{equation}}
\newcommand{\ee}{\end{equation}}

\newcommand{\pr}{\parallel}
\newcommand{\vl}{\mid}

\newcommand{\rar}{\rightarrow}

\newcommand{\lrar}{\longrightarrow}

\newcommand{\bco}{\mbox{co} \,}

\newcommand{\bsup}{\mbox{sup} \,}
\newcommand{\bmax}{\mbox{max} \,}
\newcommand{\bmin}{\mbox{min} \,}
\newcommand{\binf}{\mbox{inf} \,}

\newcommand{\meq}{\geq}

\newcommand{\al}{\alpha}

\newtheorem{thm}{Theorem}[section]

\newtheorem{rem}{Remark}[section]

\newtheorem{lem}{Lemma}[section]

\begin{document}

AMS 517.977.5

\vspace{0.5cm}

% \bc{\bf (University of
%C.Peterburgsky)} \ec \vspace{0.5cm}

\begin{center}
{\Large \bf Equivalent substitution in the control theory }
\end{center}

\vspace{1cm}

\vspace{0.2cm} \bc {\bf I.M. Proudnikov} \ec  \bc {\em
pim\underline{ }10@hotmail.com}, ph.num. +79203011393 \ec

\vspace{1cm}

In this paper a system of the differential equations with a
control is considered. We study a problem of looking for an
optimal control that gives an infimum for an optimized functional.
The system of differential equations is replaced by two systems
with the upper  and lower  envelopes of a function on the right
hand side of the initial system of the differential equations. The
optimized functional is replaced by its lower envelope. All
replacements are done in a region of attainability. The necessary
conditions of optimality are sufficient for the substituted
system. The rules for evaluation of the attainability set with the
help of positively definite  functions are given in the second
part of the paper.


\bigskip
% \vspace{1.5cm}
\normalsize \section{Introduction}

The system of the differential equations is considered,  the right
side of that includes a control $u$. The problem is formulated as
following: to find a control $u$ from a set $U$ for that a given
functional $J(\cdot)$ gets its optimum. The rules for equivalent
replacement of our system by more simple system are given for that
the functional $J(\cdot)$ gets the same optimum. Moreover, the
necessary optimal conditions turn into sufficient optimal
conditions. All replacements are done in a region of attainability
and consist in construction of the upper concave or lower convex
envelopes of the optimized functional $J(\cdot)$ and the functions
in the right side of the system. A method for construction of the
envelopes is given in Appendix.

Many specialists in different areas are interested in evaluations
of a region of attainability (see
\cite{chernousko}-\cite{kurganckivalyi2}). The author suggests a
new method for evaluation of a region of attainability. This
method is based on study of the differential equations of the
first order with $n-1$ constants where $n$  is the order of the
system of the differential equations. The solutions of the
differential equations of the first order are the projections of
the solutions of the system of the differential equations on a
chosen direction $g$ (see further). Changing the constants, the
direction $g$  and solving the differential equation of the first
order we can obtain enough precise evaluation of a region of
attainability.

Consider the following general problem of the control theory.
Suppose we have a system of the differential equations \be \dot{x}
(t) = \varphi (x (t), u (t)), \, \, \, \, \, x (0) = x_0,
\label{eqsub1} \ee $\varphi(\cdot)=(\varphi_1(\cdot)),
\varphi_2(\cdot)), \cdots, \varphi_n(\cdot))^*$, and an optimized
functional has a form \be J (u) = \int^T_0 f (x (\tau), u (\tau))
d \tau \rar \inf_{u}, \label{eqsub2} \ee where $ x (t) \in
\mathbb{R}^n $, $ u (t) $ takes values in $ U \subset
\mathbb{R}^r$, where $U$ is a convex compact set in $ \mathbb{R}^r
$, $ t \in [0, T] $. We assume that the function $ f (\cdot,
\cdot): \mathbb{R}^n \times \mathbb{R}^r \rar \mathbb{R} $ is
continuous and the function $ \varphi (\cdot, \cdot) $ is a
Lipschitz one in all arguments in totality, so that the system
(\ref{eqsub1}) satisfies the conditions of uniqueness for a
solution with the given initial values. We consider the autonomous
systems of the differential equations it does not restrict the
generality of consideration, as soon as if in the non-autonomous
case $ \varphi (\cdot) $ is Lipschitz in all arguments uniformly
in $t$, then all argumentations, made below, will be true as well.

We have to find an optimal control  $ u (\cdot) $ that is a
piecewise continuously differentiable vector-function from $ KC^1
[0, T] $ with values in $ U $. We will assume that the derivatives
$ u '(\cdot) $, where they exist, are bounded in the norm, i.e.
$$
\pr u '(t) \pr \leq C \, \, \, \forall t \in \aleph_u [0, T],
$$
uniformly in $u(\cdot)$. Here $ \aleph_u [0, T] $ is a set of
points in $ [0, T] $, where the derivatives  $ u' (\cdot) $ exist.
In this case  pointwise convergence of a sequence $ \{u_k (\cdot)
\} $ on $ [0, T] $ is equivalent to  uniform convergence of the
functions $ u_k (\cdot) $ on a set of continuity  and, of course,
is equivalent to  convergence in the metrics  $ \rho $ of the
space $KC^1 [0, T] $. The metrics of $KC^1 [0, T] $ is equal, by
definition, to the metrics of the space $ C [0, T] $, i.e.
$$
\rho (u_1 (t), u_2 (t)) \stackrel{\mathrm{def}}{=} \max_{t \in [0,
T]} \| u_1 (t) - u_2 (t) \|.
$$
The function $ u (\cdot) $ is defined on the segment $ [0, T] $.
We are looking for $ u (\cdot) $ for that a solution to the system
(\ref{eqsub1}) gives an infimum for the functional $ J (\cdot) $.

In the beginning we replace the optimization problem
(\ref{eqsub2}) by the following problem \be J (u, t) = \int^t_0 f
(x (\tau), u (\tau)) d \tau \lrar \binf_{u \in KC^1 [0, T], t \in
[0, T]}. \label{eqsub2b} \ee

Let us include into consideration all the functions $u(t) \in KC^1
[0, T]$ resulting from pointwise convergence. It is obvious that
all limit functions belong to a closed, bounded set of functions
defined on $ [0, T] $, which we denote by ${\cal KC}^1 [0, T] $.
The two functions from the set ${\cal KC}^1 [0, T] $ are
equivalent (equal) if these functions are equal on a set of full
measure. It is clear that all measurable functions on $[0,T]$
belong to ${\cal KC}^1 [0, T] $.

We will solve the above formulated optimization problem
(\ref{eqsub2b}) on the set ${\cal KC}^1 [0, T] $, i.e. \be J (u,
t) \lrar \binf_{{u \in{\cal KC}^1 [0, T]},{t \in [0, T]}}.
\label{eqsub2a} \ee

The problem is that an optimal control does not exist always. For
this reason generalized control (lower or upper semicontinuous) is
considered.

As an example, consider the following system of differential
equations
$$
\dot{x} (t) = u, \,\,\,\,  x(0)=0,
$$
and the optimized functional is defined as
$$
J (u) = \int^1_0 ((1-u^2)^2 + x^2) d \tau \rar \inf_{u \in{\cal
KC}^1 [0,1]}.
$$
This problem does not have an optimal control $ u (\cdot) $ in the
set of piecewise continuously differentiable functions on $[0,1]
$, but it has an optimizing sequence of controls $ \{ u_k
(\cdot)\} $ that are piecewise continuous functions with values $
\pm 1 $. It is easy to see that an optimizing sequence $\{ x_k
(\cdot) \}$, corresponding to the sequence of controls $ \{ u_k
(\cdot)\} $, has the limit $ x \equiv 0 $ on $ [0,1] $. The
control $u(\cdot) \equiv 0$, that corresponds to the solution
$x\equiv 0$, can not be received as the pointwise limit of $ \{
u_k (\cdot)\} $ and it's not an optimal control.

The right-hand side of the equation (\ref{eqsub1}) can be very
complex, and the exact solution of this equation can often be
found approximately  using numerical methods. Optimization of the
function $ J (\cdot) $  is also not easy if it has a complex form.
But optimization  of the lower convex envelope (LCE) of $ J
(\cdot) $ is easier. Moreover,  a global optimum point does not
disappear if we construct a lower convex envelope of our optimized
functional. In addition, the construction of the lower convex
envelope of $ J (\cdot) $, that we denote by $ \tilde{J} (\cdot)
$, turns it into a lower weakly semicontinuous function. It means
that
$$
{\underline{lim}}_{u_k \lrar u} J (u_k) \meq J (u)
$$
for any sequence $ \{u_k \} $  converging to $u$ weakly. This
requirement is important for  weak convergence of an optimizing
sequence to a solution of the problem (\ref{eqsub2}).

We propose here a method of equivalent substitution with the help
of which we can overcome these difficulties. Namely, we suggest a
replacement of the right-hand side of the equation (\ref{eqsub1})
by another function with a simpler structure. The search for
solutions of the system (\ref{eqsub1}) (numerical or not) becomes
simpler. The principle of equivalent replacement claims that
although we have another function with a simpler structure, but
the function $\tilde J (\cdot) $ attains the same infimum on the
set of piecewise continuously differentiable functions and on its
closure. At the same time the new optimized functional $ \tilde J
(\cdot) $ becomes lower semicontinuous.

This idea is different from the idea of relaxation \cite{cesary}
because we take the convex (concave) envelopes of functions on the
right hand side of our system of the differential equations for
all points from the region of attainability.

Taking into consideration the information  about the replacements
of the functions $ \varphi (\cdot) $ and $ J (\cdot) $, we can
conclude that searching for an optimal control and  optimal
trajectory becomes easier and the new optimization problem is
equivalent to the initial one in the sense of finding  an optimal
control. In this case the conditions of optimality become
necessary and sufficient.

For the first time the author used the idea about replacing a
function with its lower convex approximation for finding its
optimal points in \cite{lowapp1} -\cite{lowapp3}.

\bigskip
\section{\bf The principle of equivalent replacement}
\bigskip

Let us consider the same  system of differential equations
(\ref{eqsub1}), and the optimized functional (\ref{eqsub2a}). We
rewrite the system (\ref{eqsub1}) and (\ref{eqsub2a}) in the
following form  \be
\left \{\begin{array}{ll} \dot{x} (t) = \varphi (x (t), u (t)) , \\
\dot{y} (t) = f (x (t), u (t))
\end{array} \right. \label{eqsub3} \ee
with the initial conditions $ x (0) = x_0, \, y (0) = 0 $. The
optimization problem
 (\ref{eqsub2a}) is replaced by another optimization problem
\be y (t, u) \lrar \binf_{{u \in{\cal KC}^1 [0, T]},{t \in [0,
T]}}. \label{eqsub4} \ee Any solution of (\ref{eqsub3}) is a
solution of the integral equations \be \left \{
\begin{array}{ll} x (t) = \int^t_0 \varphi (x (\tau), u (\tau))d\tau + x_0,
\vspace{0.2cm} \\y (t) = \int_0^tf (x (\tau), u (\tau))d \tau,
\vspace{0.2cm} \\u (t) \in{{\cal KC}^1 [0, T]}, \, \, t \in [0,
T],
\end{array} \right. \label{eqsub5} \ee
where $ u (\cdot) $ is a piecewise continuously differentiable
control.

Unite all solutions of (\ref{eqsub5}) in one set $ D (t) $ for $t
\in [0,T]$

$$
D (t)  =  \{(x, y, z) \vl x= x (t) = \int^t_0 \varphi (x (\tau),
u (\tau)) d \tau + x_0,  \, \,
$$
\be
 y=y (t) = \int_0^t f (x (\tau),
u (\tau)) d \tau,{} z = u (t) \in{\cal KC}^1 [0, T] \}, \ee which
is called the set of attainability for the systems (\ref{eqsub1})
and (\ref{eqsub2a}) at time $ t \in [0, T] $.

It is easy to see that the optimization problem (\ref{eqsub4}) is
equivalent to the following optimization problem \be L (x, y, z) =
y \lrar {\binf}_{y \in{\bigcup_{t \in [0, T]}} D (t)}.
\label{eqsub6} \ee

The function $ L (x, y, z) = y $ is linear in the coordinates $
(x, y, z) $. (The function $ L (\cdot, \cdot, \cdot) $ depends
only on the coordinate  $ y $). It is well known that any linear
function reaches its maximum or minimum on boundary of any compact
set on which maximum or minimum are looked for.

Since the set of solutions of (\ref{eqsub1}) in accordance with
the assumptions is bounded on $ [0, T] $, the set of the
vector-valued functions $u(\cdot)$ is closed and bounded in ${\cal
KC}^1 [0, T] $. Then $ D (t) $ is closed and bounded for any $ t
\in [0, T] $ in the metrics $ \rho $ of the space ${\cal KC}^1 [0,
T] $. Indeed, if $u_k \xrightarrow{\rho}{ } u \in U $, then, as it
was mentioned above, there is uniform convergence of $ u_k (\cdot)
$ to $ u (\cdot) $ on $ [0, T] \backslash e $, where $e$ is a set
of any small measure. Then  convergence in measure
\cite{kolmogoroffomin} holds. The convergences
$$
\varphi (x (\tau), u_k (\tau)) \rar_k \varphi (x (\tau), u (\tau))
\, \, \, \forall \tau \in [0, T],
$$
and
$$
f (x (\tau), u_k (\tau)) \rar_k f (x (\tau), u (\tau)) \, \, \,
\forall \tau \in [0, T].
$$
follow from continuity of the functions $ \varphi (\cdot, \cdot)
$, $ f (\cdot, \cdot) $ in all variables.

Uniform convergence of the integrals
$$
\int^t_0 \varphi (x (\tau), u_k (\tau)) d \tau \rar_k \int^t_0
\varphi (x (\tau, u (\tau)) d \tau
$$
and
$$
\int^t_0 f (x (\tau), u_k (\tau)) d \tau \rar_k \int^t_0 f (x
(\tau), u (\tau)) d \tau.
$$
in $ t \in [0, T] $ and $ k $ follows from Egorov's theorem
\cite{vulich}.

Indeed, otherwise the sequences $ \{T_k \} $ and $ \{u_k (\cdot)
\} $ exist, for that and for some $ \varepsilon> 0 $, the
inequalities
$$
\mid \int^{t_{k}}_0 \varphi (x (\tau), u_k (\tau)) d \tau -
\int^{t_{k}}_0 \varphi (x (\tau, u (\tau)) d \tau \mid>
\varepsilon
$$
and
$$
\mid \int^{t_{k}}_0 f (x (\tau), u_k (\tau)) d \tau -
\int^{t_{k}}_0 f (x (\tau), u (\tau)) d \tau \mid> \varepsilon.
$$
hold. The integrals can be considered as functions of $ t $.

According to Egorov's theorem for any small $ \delta> 0 $ there is
a set $ e $  with measure $ \mu (e) <\delta $ , that  the
integrals, as the functions of $ t $, will converge uniformly in $
k $ on a set $ [0, T] \backslash \, e $. As soon as the integrals
are absolutely continuous in measure, the integrals over the set
$e$ with measure $\mu(e)<\delta $ will be arbitrarily small if $
\delta $ is also arbitrarily small. As a result, we come to the
contradiction with existence of $\varepsilon $, for which the
inequalities, written above, are true. Thus we have proved the
following theorem.

\begin{thm} The set $ D (t) $ is closed and bounded in ${\cal KC}^1 [0, T] $
for any $ t \in [0, T] $ in the metrics $ \rho $ of the space $
KC^1 [0, T] $.
\end{thm}

Consider a sequence of the functions defined on $ [0, T] $, \be
x_{k +1} (t) = \int^{t}_0 \varphi (x_{k} (\tau), u_k (\tau)) d
\tau. \label{eqsub7} \ee The sequence $ \{x_k (\cdot) \} $
converges on $ [0, T] $ uniformly in $ k $, if the sequence  $ u_k
(\cdot) $ converges in the metrics $ \rho $ to the function $ u
(\cdot) $ a.e. on $ [0, T] $. Prove this fact.

Indeed, we know from the said above that the functions $ u_k
(\cdot) $ converge to $ u (\cdot) $ uniformly on $ [0, T] $. We
replace the control $ u_k (\cdot) $ by the control $ u (\cdot) $
in (\ref{eqsub7}). The difference between the original value of
the integral (\ref{eqsub7}) and the new value of the same integral
can be evaluated in the  following way. According to the
inequality
$$
\mid \varphi (x_{k} (\tau), u_k (\tau)) - \varphi (x_{k} (\tau), u
(\tau)) \mid \leq L \mid u_k (\tau) - u (\tau) \mid \, \, \,
\forall \tau \in [0, T],
$$
where $ L $ is a Lipschitz constant of the function $ \varphi
(\cdot, \cdot) $, the mentioned above difference
 is arbitrarily small for large $ k $ as well. Indeed, we have
$$
\mid \int_0^t \varphi (x_{k} (\tau), u_k (\tau)) - \int_0^t
\varphi (x_{k} (\tau), u (\tau)) \mid \leq L \int_0^t \mid u_k
(\tau) - u (\tau) \mid \, \, \, \, \forall t \in [0, T]
$$
and the right hand side of this inequality is arbitrary small for
large $k$.

We will use the following result.
\begin{lem} \cite{tricomi}, \cite{krasnov}. The sequence
$ \{x_k (\cdot, u) \}, k = 1,2, \dots $ converges uniformly for $u
\in U $ and $ k $ to a solution of (\ref{eqsub1}).
\end{lem}

Summing up everything mentioned above, we can conclude about
uniform convergence on $ [0, T] $ of the solutions $ x_k (\cdot) $
of (\ref{eqsub1}) for $ u = u_k(\cdot) $  to a solution $ x
(\cdot) $ of the same system (\ref{eqsub1}) with the control $ u
(\cdot) $ as $ k \rar \infty $.

\begin{lem} The sequence
$ \{x_k (\cdot) \}, k = 1,2, ... , $ defined by (\ref{eqsub7}),
converges uniformly on $ [0, T] $ in $k$ to a solution $ x (\cdot)
$ of (\ref{eqsub1}). \label{lemequivsubst1}
\end{lem}

\begin{rem}
Lemma (\ref{lemequivsubst1}) is also valid for the case when $ u_k
\rar u $ in the metrics $ \rho_1 $ of the space $ L_1 [0, T] $,
i.e.
$$
\rho_1 (u_k, u) = \int^T_0 \mid u_k (\tau) - u (\tau) \mid d \tau.
$$
\end{rem}

The problem (\ref{eqsub6}) has a solution if the functional
(\ref{eqsub2}) is lower semicontinuous. It will be shown how to
make it lower semicontinuous.

If there is a solution to the problem (\ref{eqsub6}) on the set of
piecewise continuously differentiable functions $ KC^1 [0, T] $,
then we have a solution to the problem \be L (x, y, z) = y \lrar
\inf_{y \in{\bco} \bigcup_{t \in [0, T]} D (t)}, \label{eqsub8}
\ee where $ \bco $ is a symbol of taking convex hull.

We introduce a set of attainability (or an attainability set)  for
the time $ T $, which, by definition, is \be D_T = \overline{\bco}
\, \cup_{t \in [0, T]} \, D (t), \label{eqsub8a} \ee where $
\overline{\bco} $ means closed convex hull. It is easy to see that
for an arbitrary
$$
(x_k (t_k), y (t_k), u_k(t_k)) \in D (t_k)
$$
such as \[(x_k (t_k), y (t_k), u_k(t_k)) \rar_k (x (t), y (t),
u(t) )
\] and
\[t_k \in [0, T], \, t \in [0, T], \, \, \, t_k \rar_k t, \]
the inclusion \[(x (t), y (t), u(t)) \in D (t) \] will be true.
Therefore, the closure in (\ref{eqsub8a}) can be removed and
definition of the set $ D_T $ can be given as the following \be
 D_T ={\bco} \, \cup_{t \in [0, T]} \, \label{eqsub8c}
D (t). \ee
Moreover, the problems (\ref{eqsub6}) and
(\ref{eqsub8}) are equivalent that means: if one of them has a
solution, then the other one has a solution as well and these
solutions are equal to each other. In addition, since projections
of the set $ D (t), t \in [0, T] $ on the axes $ x, y $ are closed
and bounded, and, hence, compact in the corresponding
finite-dimensional spaces and $D(\cdot)$ is continuous in $ t $ as
the set-valued mapping, then $ \inf $ in (\ref{eqsub8}) can be
replaced by $ \min $ and the problem (\ref{eqsub8}) can be
rewritten in the following way \be L (x, y, z) = y \lrar \min_{y
\in  D_T}. \label{eqsub8b} \ee

But a global optimal point of the problem (\ref{eqsub8b}) will not
change if we replace the functions $\varphi_i(\cdot, \cdot), i \in
1:n,$ by their upper concave and lower convex envelopes and
$f(\cdot)$ by its lower convex envelope constructed on a set of
attainability for the time $ T $. Further we will understand under
taking the upper concave and lower convex envelopes of
$\varphi(\cdot, \cdot)$ the similar operations for all coordinates
of $\varphi(\cdot, \cdot)$.

Indeed, take two arbitrary points $ (x_1 (t), y_1 (t), u_1) $ and
$ (x_2 (t), y_2 (t), u_2), $ $ t \in [0, T] $ from the set $ D (t)
$. Consider a combination with nonnegative coefficients $ \al_1 $,
$ \al_2 $, $ \al_1 + \al_2 = 1 $. Then, the point $ (\al_1 x_1 (t)
+ \al_2 x_2 (t), \al_1 y_1 (t) + \al_2 y_2 (t), \al_1 u_1 + \al_2
u_2) $ will belong to the set $ \bco D (t) $, if we replace the
functions $ \varphi (\cdot, \cdot) $ and $ f (\cdot, \cdot) $ by
the following:
$$
\tilde{\varphi} (\al_1x_1 (\tau) + \al_2 x_2 (\tau), \al_1 u_1 +
\al_2 u_2) = \al_1 \varphi (x_1 (\tau), u_1) + \al_2 \varphi (x_2
(\tau), u_2)
$$
and
$$
\tilde{f} (\al_1x_1 (\tau) + \al_2 x_2 (\tau), \al_1 u_1 + \al_2
u_2) = \al_1 f (x_1 (\tau), u_1) + \al_2 f (x_2 (\tau), u_2).
$$
But this construction, performed for all points of the regions $ D
(t), t \in [0, T], $ it just means that we construct the lower and
upper convex envelopes of the function $ \varphi (\cdot, \cdot) $
and the lower convex envelope of the functions $f(\cdot, \cdot)$
in the attainability set for the time $ T $, i.e. in $D_T$.
Indeed, it follows from the above formula   (\ref{eqsub8c}) for
$D_T$ that the function $L(\cdot, \cdot, \cdot)$ reaches its
minimum (\ref{eqsub8b})  on some $D(\tau), \tau \in [0,T]$.
Consequently, we can construct the lower convex or upper concave
envelopes of the functions $\varphi(\cdot)$ and $f(\cdot, \cdot)$
in all region $D_T$.

Denoted by $ \tilde{J} (\cdot, \cdot) $ a new optimization
function obtained after the replacement of the function $ f
(\cdot, \cdot) $ by $ \tilde{f} (\cdot, \cdot) $ in $ D_T \, $
$$
\tilde{J} (u, t) = \int^t_0 \tilde{f} (x (\tau), u (\tau)) d \tau.
$$
It is clear that $ \tilde{J} (\cdot, \cdot) $ takes the same
optimal value in the attainability set $ D_T $, that the
functional $ (\ref{eqsub2a}) $ $ J (\cdot, \cdot) $ takes for the
system (\ref{eqsub1}).

Replace the system $ (\ref{eqsub1}) $ by the system \be \dot{x}
(t) = - \varphi (x (t), u (t)), \, \, \, \, \, x (0) =-x_0,
\label{eqsub9} \ee and the optimized functional by the functional
\be J (u, t) = \int^t_0 f (-x (\tau), u (\tau)) d \tau.
\label{eqsub10} \ee

It is easy to see that the minimum or  the maximum of the
functional $ J (\cdot, \cdot) $ did not change. Hence, the
problems $ (\ref{eqsub1}) $, $ (\ref{eqsub2a}) $ and $
(\ref{eqsub9}) $, $ (\ref{eqsub10}) $ are replaceable. So
"convexification" $ \, $ of the function $ \varphi (\cdot, \cdot)
$, in contrast to the procedure of "convexification" of the
function $ f (\cdot, \cdot) $ should be as following:
\begin{enumerate}
    \item Construction of the lower convex envelope (LCE) of the function
    $ \varphi_i (\cdot, \cdot) $ in the variables $ (x, u) $ for each
    $i \in 1:n$ from the
    attainability set for the time $T$, i.e. $ D_T $,
    which we denote by $ \tilde{\varphi}_{1i} (\cdot, \cdot) $. LCE
    of $ \varphi_i (\cdot, \cdot) $ is the biggest convex function
    that does not exceed $ \varphi_i (\cdot, \cdot) $ in $D_T$.
    \item Construction of the upper concave envelope (UCE) of the  function
     $ \varphi_i (\cdot, \cdot) $ (or, equivalently, we construct the
     lower convex envelope for the function
     $ - \varphi_i (\cdot, \cdot) $ and after that take minus of this
     function) in the
    variables $ (x, u) $ for each $i \in 1:n$  from the
    attainability set for the time $T$, i.e. $ D_T $, which we denote
    by $ \tilde{\varphi}_{2i} (\cdot, \cdot) $. UCE of $\varphi_i(\cdot, \cdot)$
    is the smallest concave function that is not less $\varphi_i(\cdot,\cdot)$
    in $D_T$.
    \item Let us replace the system $ (\ref{eqsub1}) $, $ (\ref{eqsub2a}) $
    by two systems of the equations:
\be \dot{x} (t) = \tilde{\varphi}_1 (x (t), u (t)), \, \, \, \, \,
x (0) = x_0, \, \, u (\cdot) \in{\cal KC}^1 [0, T] \label{eqsub11}
\ee with the optimization function $ \tilde{J} (u, t) $ and \be
\dot{x} (t) = \tilde{\varphi}_2 (x (t), u (t)), \, \, \, \, \, x
(0) = x_0, \, \, u (\cdot) \in{\cal KC}^1 [0, T] \label{eqsub12}
\ee with the same optimization function $ \tilde{J} (u,t) $;
   \item Let us find among the solutions of $ (\ref{eqsub11}) $ and
   $ (\ref{eqsub12}) $
    such that gives the smallest value of the
   functional $ \tilde{J} (u, t) $ in $ D_T $.
\end{enumerate}
We obtain the following result.
\begin{thm}
There are the solutions among the solutions of $ (\ref{eqsub11}) $
and $ (\ref{eqsub12}) $  such that deliver a minimum (maximum) in
$u(\cdot) \in {\cal KC}^1 [0, T]$ and $t \in [0,T]$   for the
functional
$$
\tilde{J} (u, t) = \int^t_0 \tilde{f} (x (\tau), u (\tau)) d \tau,
$$
that coincides with an infimum (supremum) of the functional $ J
(u, t) $ (see (\ref{eqsub2a})). Moreover, necessary conditions for
the minimum (maximum) are also sufficient conditions.
\end{thm}

\begin{rem}
The set $ D (T) $ is not necessarily compact, although its
projections on the axis $ x, y $ are compact . That's why we are
able to go  to the problem
$$
L (x, y, z) = y \lrar \min_{y \in D_T},
$$
if the problem (\ref{eqsub6}) has a solution. The last one
coincides with the formulation of Mazur's theorem. It asserts that
in any weakly convergent sequence $ \{u_k (\cdot) \} \in L_p ([0,
T]) $, $ u_k (\cdot) \longrightarrow u (\cdot), $ a subsequence
can be chosen for each $ k $  convex hull of which is almost
everywhere on $ [0, T] $ converges as $ k \rar \infty $ to some $
u (\cdot) \in L_p [0, T] $. In our case, there exists a sequence $
\{u_k (\cdot) \} \in {\cal KC}^1 [0, T] $, the convex hull of
which will converge to an optimal control $ u (\cdot) \in {\cal
KC}^1 [0, T] $. The sequence of the solutions $ \{x_k (\cdot) \}
$, corresponding to the controls  $ u_k (\cdot)$,  will converge
to an optimal solution $ x (\cdot) $, corresponding to the control
$ u (\cdot) $, provided that the solutions  have been calculated
to the problems with the modified right-hand side.

\end{rem}

\begin{rem}

The rules for construction of LCE and UCE are given in Appendix.

\end{rem}

\begin{rem}

In many cases we have to construct only LCE or UCE  for the
function $ \varphi(\cdot)$.

\end{rem}

Return back to the initial problem (\ref{eqsub2}) with the fixed
time $T$. Consider a set
$$
 D(t)  =  \{ (x,y,z) \vl x=x(t)= \int^t_0
\varphi(x(\tau),u(\tau))d\tau+ x_0,
$$
$$
y= y(T)= \int_0^T f(x(\tau),u(\tau))d\tau, {}
      z=u(t) \in {\cal KC}^1 [0,T] \},
%\end{eqnarray}
$$
that is called the set of attainability of the system
(\ref{eqsub1}),(\ref{eqsub2}) at time $t$.

Let us introduce a set of attainability for the time $T$ for the
system (\ref{eqsub1}),(\ref{eqsub2}) that is by definition
\be
D_T= \overline{\bco} \, \cup_{t \in [0,T]} \, D(t).
\label{eqsub12a} \ee

As above it is possible to prove that we can remove the closure in
(\ref{eqsub12a}) and write
$$
D_T= {\bco} \, \cup_{t \in [0,T]} \, D(t).
$$
The optimization problem can be reformulated in the form \be
L(x,y,z)=y \lrar \inf(\sup)_{(x, y, z) \in D_T}. \label{eqsub13a}
\ee

The problems (\ref{eqsub1}), (\ref{eqsub2}) and (\ref{eqsub13a})
are equivalent which means if one has a solution, then another one
has a solution and these solutions are the same. Moreover, as soon
as the projections of the sets $D(t), t \in [0,T],$ on the axes
$x,y$ are closed, bounded and continuous as a set valued mappings,
then we can write instead of $\binf, \bsup$ $\bmin, \bmax$ if a
solution of (\ref{eqsub13a}) exists.

We come to the following result.
\begin{thm}
There are some solutions among the solutions of $ (\ref{eqsub11})
$ and $ (\ref{eqsub12}) $  such that deliver a minimum (maximum)
in $u(\cdot) \in {\cal KC}^1 [0, T]$ and $t \in [0,T]$ to the
functional
$$
\tilde{J} (u) = \int^T_0 \tilde{f} (x (\tau), u (\tau)) d \tau,
$$
that coincides with an infimum (supremum) of the functional $ J
(u) $ (see (\ref{eqsub2})) where $ \tilde{f} (\cdot, \cdot) $.is
LCE of the function $ {f} (\cdot, \cdot) $ Moreover, the necessary
conditions for minimum (maximum) are also the sufficient
conditions.
\end{thm}

%It is clear that the attainability set $D_T$ will not change if we
%instead of the systems (\ref{eqsub11}) and (\ref{eqsub12})
%consider the systems
%$$
%\dot{x} (t) = \sum_{i=1}^n \al_i \tilde{\varphi}_{1i} (x (t), u
%(t)), \, \, \, \, \, x (0) = x_0, \, \, u (\cdot) \in{\cal KC}^1
%[0, T]
%$$
%and
%$$
%\dot{x} (t) = \sum_{i=1}^n \al_i \tilde{\varphi}_{2i} (x (t), u
%(t)), \, \, \, \, \, x (0) = x_0, \, \, u (\cdot) \in{\cal KC}^1
%[0, T]
%$$
%for different $\al_i \meq 0$, $\sum_{i=1}^n \al_i =1$,
%$\tilde{\varphi}_{1i}(\cdot), \tilde{\varphi}_{2i}(\cdot)$ are the
%coordinates of the vector-functions $\tilde{\varphi}_{1}(\cdot),
%\tilde{\varphi}_{2}(\cdot)$ correspondingly.

Consider some examples. It is clear that an equivalent replacement
of one system by another can be applied to a differential system
without control $u$.

Example 1. Consider the differential equation
$$
\dot{x} (t) = \varphi (x (t)) = \left \{
\begin{aligned}
(X-1)^2, \mbox{if $ x \meq 0 $} \\
(X +1)^2, \mbox{if $ x <0 $} \\
\end{aligned} \right.
$$
with the initial condition $ x (0) = 0 $. The optimized functional
is given by
$$
f (x (t)) = x^2 (t) \rar \min \, \, \, \mbox{for $ t \in (-
\infty, + \infty) $},
$$
The general solution of the differential equations for $ x \meq 0
$ has the form
$$
x (t) = - \frac{1}{t + c} +1,
$$
which tends to $ 1 $ as $ t \rar \infty $. The general solution of
the  differential equation for $ x <0 $ is given by
$$
x (t) = - \frac{1}{t-c} -1,
$$
which tends to $ -1 $ as $ t \rar \infty $. In order to meet the
initial condition we have to put $ c =  1$. The projection of the
attainability set $ D_{(- \infty, + \infty)} $ on $ OX $ axis is
the interval $ (-1, +1).$

It is clear that the function $ f (\cdot) $ takes its minimum at $
x = 0 $. But we get the same solution  if instead of the function
$ \varphi (\cdot) $ we take its lower convex envelope, namely, the
function
$$
\tilde{\varphi} (x) =
\begin{cases}
(X-1)^2, \mbox{if $ x \meq 1 $} \\
0, \mbox{if $ -1 \leq x \leq 1 $} \\
(X +1)^2, \mbox{if $ x <-1 $}. \\
\end{cases}
$$

Example 2. The same example, but
$$
{J} (x, u) = \int^t_0 x^2 (\tau)) d \tau,
$$
which we minimize for $ t \in [0,1] $. The equation of the
solution is
$$
x (t) = - \frac{1}{t +1} + 1.
$$
Here the replacement of the function $ \varphi (\cdot) $ by the
function $ \tilde{\varphi} (\cdot) $ on the whole line is not
correct, since the projection of the attainability set $ D_{1} $
on $ OX $ axis is the interval $ [0, 1/2] $.

Example 3. Let us give the differential equation
$$
\dot{x} (t) = x^2
$$
with the initial condition $ x (0) = 1 $. The general solution has
the form
$$
x (t) = - \frac{1}{t + c},
$$
a solution, satisfying the initial condition, is
$$
x (t) = - \frac{1}{t-1}.
$$
The optimized functional has the form
$$
{J} (x, u) = \int^t_0 (-x^2 (\tau)) d \tau \longrightarrow_t \, \,
\binf
$$
for $ t \in [0,1] $.

 It is easy to compute its optimal value
$$
{J} (x, u) = \int^1_0 (-x^2 (\tau)) d \tau = \int^1_0 (- \dot{x}
(\tau)) d \tau = x (0)-x (1) = - \infty.
$$

In this case, the projection of the attainability set $ D_{1} $ on
$ OX $ axis is the set $ (- \infty, 0) \cup [1, + \infty) $. It is
easy to see that the lower convex envelope of the functional $ J
(\cdot) $ on $ D_{1} $, which we denote by $ \tilde{J} (\cdot) $,
takes the same infimum value. It is also true for the functional
$$
{J} (x, u) = \int^{\infty}_1 (-x^2 (\tau)) d \tau =
\int^{\infty}_1 (- \dot{x} (\tau)) d \tau = x (1)-x (+ \infty) = -
\infty.
$$

Example 4. Let us consider the following problem
$$
\dot {x}(t)= x sin(1/x) + u, \,\,\, x(0)=0.
$$
The optimized functional is
$$
J(u)=\int_0^{\infty} \vl u(\tau)-x(\tau) \vl d \tau \rar \binf_u.
$$
We will get the following system after construction of the lower
and upper envelopes
$$
\dot {x}(t)= x  + u, \,\,\, x(0)=0
$$
and
$$
\dot {x}(t)= - x  + u, \,\,\, x(0)=0.
$$
The optimal solution exists among their solutions $x(t)\equiv
u(t)= 0$.

Example 5. Let us consider the differential equation
$$
\dot {x}(t)= x^2 - u^2
$$
with the initial condition $x(0)=0$. We are considering piecewise
continuously differentiable functions $u(\cdot), |u(\cdot)| \leq
1, $ on the segment $[0,1]$ for such that delivers minimum to the
functional
$$
J(u)=\int_0^1 x^2(\tau) d\tau .
$$
The solution when $u(\cdot)$ is constant on the segment $[0,1]$ is
given by the form
$$
x(t)=\frac{u(1 - e^{2u(t+c)})}{1 + e^{2u(t+c)}}.
$$
Here the constant $c$ is defined by the initial conditions. We can
see from here that if $u(\cdot)$ is not constant on $[0,1]$, then
$| x(\cdot) | \leq | u(\cdot) |$ for any initial conditions. It
means that a curve $x(t), t \in [0,1],$ will be in a set bounded
by the lines $x=\pm u$ on the plane $XOU$, where $| x(\cdot) |
\leq | u(\cdot) |$. The set of attainability $D_T, T=1,$ will
belong to the same set. UCE of the function $\varphi(x.u) = x^2 -
u^2$ in  $D_T$ is a function the graph of which  goes through the
point $(0,0,0)$. Therefore, if we solve the differential equations
with the right sides $\tilde{\varphi}_1(\cdot.\cdot)$ and
$\tilde{\varphi}_2(\cdot.\cdot)$, then among the solutions there
are  such that deliver the minimum $0$ to the functional
$J(\cdot)$ i.e. the formulated theorem is true.

\bigskip
\section{\bf An evaluation of the attainability set}
\bigskip

Let us have a system of differential equations \be \dot{x} (t) =
\varphi (x, u), \, \, \, x \in \mathbb{R}^n, \, \, t \in [0, T],
\, \, u (t) \in U \subset \mathbb{R}^r \label{eqsub13} \ee with
the initial condition $ x (0) = 0 $, where $ \varphi (\cdot,
\cdot) $ is Lipschitz in the variables $ x, u $, $ U  $ is a
convex compact set in $ \mathbb{R}^r $. The problem  is to
estimate the attainability set. By definition, the area of
attainability for the time $ T $ is the set
$$
D_T = \overline{\bco} \, \cup_{t \in [0, T]} \, D (t),
$$
where
$$
D (t) = \{x \in \mathbb{R}^n \vl x = x (t) = \int^t_0 \varphi (x
(\tau), u (\tau)) d \tau, \, \, u (\tau) \in U, u (\cdot) \in
{\cal KC}^1 [0, T] \}.
$$
The choice of the initial position and the initial time of zero is
not a loss of generality.

Take an arbitrary   positively definite function $ V (x) $ (see
\cite{zub1}), satisfying the condition
$$
m_1 \| x \|^2 \leq V (x) \leq m_2 \| x \|^2.
$$
Let
$$
\varphi (x, u, t) = \varphi_1 (x, u, t) + \varphi_2 (x, u, t)
$$
and $v(\cdot,\cdot): [0,T] \times U  \rar \mathbb{R}^n$ is a
piecewise continuous vector-function.

Consider the systems of differential equations \be \dot{x} (t) =
\tilde{\varphi}_1 (x, u) = \varphi_1 (x, u)+v(u,t) , \, \, t \in
[0, T], \label{eqsub14} \ee and \be \dot{x} (t) =
\tilde{\varphi}_2 (x, u) = \varphi_2 (x, u)-v(u,t). \, \, t \in
[0, T], \label{eqsub15} \ee

We denote  by \be D_T^{(i)} = \overline{\bco} \, \cup_{t \in [0,
T]} \, D_i (t), \, i=1,2, \label{eqsub16} \ee the attainability
sets for the systems (\ref{eqsub14}), (\ref{eqsub15}), where
$$
D_i (t) = \{x \in \mathbb{R}^n \vl x = x_i (t) = \int^t_0
\tilde{\varphi}_i (x (\tau), u (\tau)) d \tau, \, \, u (\tau) \in
U,
$$
\be
 u (\cdot) \in{\cal KC}^1 [0, T] \}, i = 1,2.
\label{eqsub17} \ee Let the estimates of the attainability sets
for the time $ T $ be given respectively by the inequalities
$$
0 \leq V (x) \leq c_{1}, \, \, \, 0 \leq V (x) \leq c_{2}.
$$
We get the estimation for the attainability set of the system
(\ref{eqsub13}) for the time $ T $. We show that  the
attainability set $ D_T $ for this system satisfies the inclusion
$$
D_T \subset D_T^{(1)} + D_T^{(2)}.
$$
Indeed, by definition, the set $ D_T^{(1)} + D_T^{(2)} $ will
consist of the points on the curves the tangents to which are the
sum of the tangents to the curves consisting of the points of the
sets $ D_{1} (t_1) $ and $ D_{2} (t_2) $ for all $ t_1, t_2 \in
[0, T] $. It is clear, that for some vector-function
$v(\cdot,\cdot)$ the resulting set will include $ D_T $, which
consists of the points on the curves the tangents to which are the
sum of the tangents to the curves consisting of the points of the
sets $ D_{1} (t) $ and $ D_{2} (t) $ for all $ t \in [0, T] $.

As a result, the following theorem is proved.
\begin{thm}
For the attainability set $ D_T $ of (\ref{eqsub13}) the inclusion
$$
D_T \subset D_T^{(1)} + D_T^{(2)},
$$
is true for some vector-function $v(\cdot,\cdot)$, where $
D_T^{(i)} $, $ i = 1,2, $ are given by (\ref{eqsub16}),
(\ref{eqsub17}) .
\end{thm}

The following lemma follows from here.

\begin{lem}
The function $ V (\cdot) $ satisfies the inequality
$$
0 \leq V (x) \leq c_1 + c_2
$$
in the attainability set $ D_T $ of the system (\ref{eqsub13}).
\end{lem}

Now consider two differential systems  with the right sides
$\varphi(\cdot,\cdot)$, $\varphi_1(\cdot,\cdot)$ and zero initial
conditions at the time  equaled to zero.

Let us suppose that the vectors $\varphi(\cdot,\cdot)$ ,
$\varphi_1(\cdot,\cdot)$ are collinear and also the inequality

\be \| \varphi (x, u) \| \leq k_2 \| \varphi_1 (x, u) \|
\label{eqsub19} \ee is true for all $ x, u $. We assume that we
know the attainability set $ D^{(1)}_T $ for the time $ T $ of the
system with the right hand side $\varphi_1(\cdot,\cdot)$. The
problem is to obtain some estimates of the attainability set of
the system (\ref{eqsub13}). The arguments will be carried out as
previously, considering the trajectories of the corresponding
systems.

Any vector in the set $ D_T $ for some $ t \in [0, T] $ and $ u
\in {\cal KC}^1[0,T] $ is
$$
x = x (t) = \int^t_0 \varphi (x (\tau), u (\tau)) d \tau.
$$
Consequently,
$$
\| x (t) \| = \| \int^t_0 \varphi (x (\tau), u (\tau)) d \tau \|
\leq \int^t_0 \| \varphi (x (\tau), u (\tau)) \| d \tau \leq
$$
$$
\leq k_2 \int^t_0 \| \varphi_1 (x (\tau), u (\tau)) \| d \tau
\subset k_2 D_1 (t).
$$
Since the previous inclusion holds for any $ t \in [0, T] $, it
follows that
$$
D_T \subset \cup_{t = 0}^T k_2 D_1 (t) = k_2 \cup_{t = 0}^T D_1
(t)
 \subset k_2 \overline{co} \cup_{t = 0}^T D_1 (t) = k_2 D_T^{(1)}.
$$
The  following theorem is proved.
\begin{thm}

For the systems of the differential equations with the right sides
$\varphi(\cdot,\cdot)$, $\varphi_1(\cdot,\cdot)$ and with the
attainability sets $ D_T $ and $ D_T^{(1)} $ respectively, for
which the inequality (\ref{eqsub19}) holds, the inclusion
$$
D_T \subset k_2 D_T^{(1)}
$$
is true.
\end{thm}
From here we can easily obtain the following conclusion.
\begin{lem}
In the attainability set  $ D_T $ of the system (\ref{eqsub13})
the function $ V (\cdot) $ satisfies the inequality
$$
0 \leq V (x) \leq k_2 c_1,
$$
where the constant $ c_1 $ limits the top value of the function $
V (\cdot) $ in the attainability set $ D_T^{(1)} $ of the system
(\ref{eqsub15}).
\end{lem}

Let us give a general method for evaluation of $D_T$ of the system
(\ref{eqsub13}). This method does not require any additional
information for the system (\ref{eqsub13}).

As is known, a convex set can be given by its extreme points.
There are no problems if there is a finite number of such points.
But very often these points are unknown or their number is
infinite. We can reconstruct a convex set if we know its
projections on different directions. If we project any trajectory,
then we project not only the points but the tangents constructed
at these points. It means that we have to consider the following
system for any direction $g \in \mathbb{R}^n, \| g \| =1,$
$$ (\dot{x},g)(t)=(\varphi((x,g),u),g), \,\,\, x(t) \in \mathbb{R}^n,
\,\, t \in [0,T], \,\, u(t) \in U \subset \mathbb{R}^r. $$

As a result, we have \be \dot{\theta}(t)=(\varphi(x(t),u),g),
\,\,\, \theta(t) \in \mathbb{R}^1, \,\, t \in [0,T], \,\, u(t) \in
U \subset \mathbb{R}^r, \label{eqsub19a} \ee where $\theta$ is the
scalar production $(x,g)$.

We can do an orthogonal transformation that the direction of the
vector $g$ was the first coordinate axis $x_1$. Then the
projection of the velocity vector $(\dot{x}_1,
\dot{x}_2,\dot{x}_3, \dots, \dot{x}_n)$ on the line $x_1$  will be
equal to $(\dot{x}_1, 0,0, \dots, 0)$. It means that we have to
substitute $x_2=c_2, x_3=c_3, \dots, x_n=c_n$ into the first
equation and to solve the first order differential equation for
different values of the constants $c_2, c_3, \dots, c_n$ that are
corresponding  coordinates of the start point.

We can make the following conclusion: {\em convex hull of the
attainability sets of the equation (\ref{eqsub19a}) for different
$g \in S_1^{n-1}(0)$ will include the attainability set of the
equation (\ref{eqsub13})}. It is possible to do, because
calculation methods are developed very well for the first order
differential equations.

\bigskip
% \vspace{1.5cm}
\normalsize \section{Conclusion}

The obtained results allow us to pass from local to  global
optimization  problem. To implement this it is required to
construct the lower convex  and upper concave envelops of the
function written on the right hand side of the differential system
(\ref{eqsub1}). We also construct the lower convex envelops for
the optimized functional. All constructions are done in the
attainability set for the time $T$.

A method for  estimation of the attainability set with the help of
the positively definite  functions (Lyapunov functions) is
suggested. The proposed method is based on the decomposition of
the function, stayed on the right hand side of the system of the
differential equations, into the components the sets of
attainability of which are already known. It makes it different
from the paper \cite{kostousova}, where the linear systems are
considered.

It is suggested to find projections of $D_T$ onto any direction $g
\in \mathbb{R}^n$. For this reason we have to find  projections of
the trajectories of the differential system and the tangents to
them to the direction $g$. We come across the problem of
definition of a set of attainability for the differential equation
of the first order.

The proposed transformation method of the systems is especially
useful when it is difficult to get a solution of  differential
equations in an explicit form, but while using approximate methods
only. In addition, the sufficient conditions of optimality for an
optimal control are obtained according to the proposed method.

%\bigskip
\vspace{1cm}
%\normalsize \section{Appendix}

{\bf \large APPENDIX}

\vspace{0.5cm}

We will prove a theorem giving a rule for construction of LCE and
UCE.

Let $f(\cdot):\mathbb{R}^n \rar \mathbb{R}- $ be continuous
function on a convex compact set $D$. It is required to construct
LCE and UCE in $D$. Consider a function
$\varphi_p(\cdot):\mathbb{R}^n \rar \mathbb{R}$
$$
\varphi_p(x)=\frac{1}{\mu(D)} \int_D f(x+y)p(y) dy,
$$
where $p(\cdot)$ is a distribution function satisfying the
following equalities \be p(y) \meq 0 \;\;\;\;\; \forall y \in D,
\;\;\;\;\; \frac{1}{\mu(D)} \int_D p(y) dy = 1, \;\;\;\; \int_D y
p(y) dy =0. \label{eqsub19} \ee

We will consider the functions $\varphi_p(\cdot)$ for different
distributions $p(\cdot)$.

\begin{thm}
The functions
$$
\overline{\varphi}(x)=\bsup_{p(\cdot)} \varphi_p(x), \;\;\;\;\;
\underline{\varphi}(x)=\binf_{p(\cdot)} \varphi_p(x)
$$
are UCE and LCE of $f(\cdot)$ on $D$ correspondingly.
\end{thm}
{\bf Proof.} Without loss of generality we will consider that
$f(y) \meq 0 $ for all $ y \in D$. Divide $D$ into subsets $\Delta
D_i, i \in 1:N,$ $D=\cup_i \Delta D_i$, $\mu(D)=\sum_i \mu(\Delta
D_i)$. We can approximate the function $\varphi_p(\cdot)$ with any
precision by the integral sums \be \sum_{i=1}^N f(x+y_i) \alpha_i
\beta_i, \label{eqsub20} \ee where
$$
\alpha_i=\frac{\mu(\Delta D_i)}{\mu(D)},\,\,  \beta_i=p(y_i),
\,\,\, y_i \in \Delta D_i.
$$
It follows from (\ref{eqsub19}) that \be \sum_{i=1}^N \alpha_i
\beta_i \simeq 1, \;\;\;\; \sum_{i=1}^N y_i \alpha_i \beta_i
\simeq0. \label{eqsub21} \ee The sign $\simeq$ means that the
values  on the left hand side from this sign can be  close to the
values on the right hand side with any precision depending on $N$.
The expression (\ref{eqsub20}) means that we take a convex hull of
$N$ vectors $(x+y_1), (x+y_2), \cdots , (x+y_N)$ with coefficients
$(\alpha_1 \beta_1, \alpha_2 \beta_2, \cdots, \alpha_N \beta_N)$,
i.e. we calculate a vector
$$
\bar{x}=\sum_{i=1}^N (x+y_i) \alpha_i \beta_i \simeq
x+\sum_{i=1}^N y_i\alpha_i\beta_i \simeq x
$$
and define a value of the function $\varphi_p(\cdot)$ at this
point equaled to
$$
\sum_{i=1}^N f(x+y_i) \alpha_i \beta_i .
$$
Changing the points $x+y_i \in D$ and the coefficients $\{
\alpha_i \beta_i \},  i \in 1:N,$ satisfying (\ref{eqsub21}), we
define in such way the  functions $\varphi_p(\cdot)$ with
different values at $x$.

Let us prove that the function
$$
\underline{\varphi}(x)=\binf_{p(\cdot)} \varphi_p(x)
$$
is LCE. As soon as $\binf$ is taken for all distributions
$p(\cdot)$, then the inequality $ \underline{\varphi}(x) \leq
f(x)$ is true for all $x \in D$. The function $\varphi_p(\cdot)$
 can be approached by the sums
(\ref{eqsub20}) for any distribution $p(\cdot)$  under conditions
on the coefficients (\ref{eqsub21}). It follows from here that
$\varphi_p(\cdot)$ can not be smaller than LCE of $f(\cdot)$. The
operation $\binf$ keeps this quality. Consequently,
$\underline{\varphi}(\cdot)$ is LCE of $f(\cdot)$. We can prove in
the same way that $\overline{\varphi}(\cdot)$ is UCE of
$f(\cdot)$. The Theorem is proved. $\Box$

The construction of LCE can be done using Fenchel-Morrey's theorem
\cite{alextihomfomin}. According to it LCE is equal to the second
conjugate function $f^{**}(\cdot)$. Construction of
$f^{**}(\cdot)$ is not easy. To find a value of LCE at one point
we have to solve two difficult optimization problems, namely,
$$
f^*(p)=\bsup_{x \in D} \{ (p,x) - f(x) \}
$$
and
$$
f^{**}(x)=\bsup_{p \in D^*} \{ (p,x) - f^*(p) \}.
$$

\end{document}